\documentclass[12pt]{article}

\usepackage{amsmath,amssymb,amsthm,amscd,a4wide,upref}

\usepackage[enableskew]{youngtab}

\usepackage{hyperref}
\hypersetup{colorlinks,citecolor=blue,filecolor=black,linkcolor=blue,urlcolor=blue}
\usepackage{enumerate}
\usepackage{mathtools}
\usepackage{enumitem}

\usepackage{lmodern}     
\usepackage[T1]{fontenc}

\begin{document}

\newcommand{\End}{{\rm{End}\ts}}
\newcommand{\Hom}{{\rm{Hom}}}
\newcommand{\Mat}{{\rm{Mat}}}
\newcommand{\ad}{{\rm{ad}\ts}}
\newcommand{\ch}{{\rm{ch}\ts}}
\newcommand{\chara}{{\rm{char}\ts}}
\newcommand{\ind}{{\rm{ind}\ts}}
\newcommand{\diag}{ {\rm diag}}
\newcommand{\pr}{^{\tss\prime}}
\newcommand{\non}{\nonumber}
\newcommand{\wt}{\widetilde}
\newcommand{\wh}{\widehat}
\newcommand{\ot}{\otimes}
\newcommand{\la}{\lambda}
\newcommand{\ls}{\ts\lambda\ts}
\newcommand{\La}{\Lambda}
\newcommand{\De}{\Delta}
\newcommand{\al}{\alpha}
\newcommand{\be}{\beta}
\newcommand{\ga}{\gamma}
\newcommand{\Ga}{\Gamma}
\newcommand{\ep}{\epsilon}
\newcommand{\ka}{\kappa}
\newcommand{\vk}{\varkappa}
\newcommand{\vt}{\vartheta}
\newcommand{\si}{\sigma}
\newcommand{\vp}{\varphi}
\newcommand{\de}{\delta}
\newcommand{\ze}{\zeta}
\newcommand{\om}{\omega}
\newcommand{\ee}{\epsilon^{}}
\newcommand{\su}{s^{}}
\newcommand{\hra}{\hookrightarrow}
\newcommand{\ve}{\varepsilon}
\newcommand{\ts}{\,}
\newcommand{\vac}{\mathbf{1}}
\newcommand{\di}{\partial}
\newcommand{\qin}{q^{-1}}
\newcommand{\tss}{\hspace{1pt}}
\newcommand{\Sr}{ {\rm S}}
\newcommand{\U}{ {\rm U}}
\newcommand{\BL}{ {\overline L}}
\newcommand{\BE}{ {\overline E}}
\newcommand{\BP}{ {\overline P}}
\newcommand{\AAb}{\mathbb{A}\tss}
\newcommand{\CC}{\mathbb{C}\tss}
\newcommand{\KK}{\mathbb{K}\tss}
\newcommand{\QQ}{\mathbb{Q}\tss}
\newcommand{\SSb}{\mathbb{S}\tss}
\newcommand{\ZZ}{\mathbb{Z}\tss}
\newcommand{\X}{ {\rm X}}
\newcommand{\Y}{ {\rm Y}}
\newcommand{\Z}{{\rm Z}}
\newcommand{\Ac}{\mathcal{A}}
\newcommand{\Lc}{\mathcal{L}}
\newcommand{\Mc}{\mathcal{M}}
\newcommand{\Pc}{\mathcal{P}}
\newcommand{\Qc}{\mathcal{Q}}
\newcommand{\Tc}{\mathcal{T}}
\newcommand{\Sc}{\mathcal{S}}
\newcommand{\Bc}{\mathcal{B}}
\newcommand{\Ec}{\mathcal{E}}
\newcommand{\Fc}{\mathcal{F}}
\newcommand{\Hc}{\mathcal{H}}
\newcommand{\Uc}{\mathcal{U}}
\newcommand{\Vc}{\mathcal{V}}
\newcommand{\Wc}{\mathcal{W}}
\newcommand{\Yc}{\mathcal{Y}}
\newcommand{\Ar}{{\rm A}}
\newcommand{\Br}{{\rm B}}
\newcommand{\Ir}{{\rm I}}
\newcommand{\Fr}{{\rm F}}
\newcommand{\Jr}{{\rm J}}
\newcommand{\Or}{{\rm O}}
\newcommand{\GL}{{\rm GL}}
\newcommand{\Spr}{{\rm Sp}}
\newcommand{\Rr}{{\rm R}}
\newcommand{\Tr}{{\rm T}}
\newcommand{\Zr}{{\rm Z}}
\newcommand{\gl}{\mathfrak{gl}}
\newcommand{\middd}{{\rm mid}}
\newcommand{\ev}{{\rm ev}}
\newcommand{\Pf}{{\rm Pf}}
\newcommand{\Norm}{{\rm Norm\tss}}
\newcommand{\oa}{\mathfrak{o}}
\newcommand{\spa}{\mathfrak{sp}}
\newcommand{\osp}{\mathfrak{osp}}
\newcommand{\g}{\mathfrak{g}}
\newcommand{\h}{\mathfrak h}
\newcommand{\n}{\mathfrak n}
\newcommand{\z}{\mathfrak{z}}
\newcommand{\Zgot}{\mathfrak{Z}}
\newcommand{\p}{\mathfrak{p}}
\newcommand{\sll}{\mathfrak{sl}}
\newcommand{\agot}{\mathfrak{a}}
\newcommand{\qdet}{ {\rm qdet}\ts}
\newcommand{\Ber}{ {\rm Ber}\ts}
\newcommand{\HC}{ {\mathcal HC}}
\newcommand{\cdet}{ {\rm cdet}}
\newcommand{\tr}{ {\rm tr}}
\newcommand{\gr}{ {\rm gr}}
\newcommand{\str}{ {\rm str}}
\newcommand{\loc}{{\rm loc}}
\newcommand{\Gr}{{\rm G}}
\newcommand{\sgn}{ {\rm sgn}\ts}
\newcommand{\ba}{\bar{a}}
\newcommand{\bb}{\bar{b}}
\newcommand{\bi}{\bar{\imath}}
\newcommand{\bj}{\bar{\jmath}}
\newcommand{\bk}{\bar{k}}
\newcommand{\bl}{\bar{l}}
\newcommand{\hb}{\mathbf{h}}
\newcommand{\Sym}{\mathfrak S}
\newcommand{\fand}{\quad\text{and}\quad}
\newcommand{\Fand}{\qquad\text{and}\qquad}
\newcommand{\For}{\qquad\text{or}\qquad}
\newcommand{\OR}{\qquad\text{or}\qquad}
\newcommand{\emp}{\mbox{}}
\newcommand{\atopn}[2]{\genfrac{}{}{0pt}{}{#1}{#2}}

\renewcommand{\theequation}{\arabic{section}.\arabic{equation}}

\newtheorem{thm}{Theorem}[section]
\newtheorem{lem}[thm]{Lemma}
\newtheorem{prop}[thm]{Proposition}
\newtheorem{cor}[thm]{Corollary}
\newtheorem{conj}[thm]{Conjecture}
\newtheorem*{thm-ref}{Theorem}
\newtheorem*{mthm}{Main Theorem}
\newtheorem*{mthma}{Theorem A}
\newtheorem*{mthmb}{Theorem B}

\theoremstyle{definition}
\newtheorem{defin}[thm]{Definition}

\theoremstyle{remark}
\newtheorem{remark}[thm]{Remark}
\newtheorem{example}[thm]{Example}

\newcommand{\bth}{\begin{thm}}
\renewcommand{\eth}{\end{thm}}
\newcommand{\bpr}{\begin{prop}}
\newcommand{\epr}{\end{prop}}
\newcommand{\ble}{\begin{lem}}
\newcommand{\ele}{\end{lem}}
\newcommand{\bco}{\begin{cor}}
\newcommand{\eco}{\end{cor}}
\newcommand{\bde}{\begin{defin}}
\newcommand{\ede}{\end{defin}}
\newcommand{\bex}{\begin{example}}
\newcommand{\eex}{\end{example}}
\newcommand{\bre}{\begin{remark}}
\newcommand{\ere}{\end{remark}}
\newcommand{\bcj}{\begin{conj}}
\newcommand{\ecj}{\end{conj}}

\newcommand{\bal}{\begin{aligned}}
\newcommand{\eal}{\end{aligned}}
\newcommand{\beq}{\begin{equation}}
\newcommand{\eeq}{\end{equation}}
\newcommand{\ben}{\begin{equation*}}
\newcommand{\een}{\end{equation*}}

\newcommand{\bpf}{\begin{proof}}
\newcommand{\epf}{\end{proof}}

\def\beql#1{\begin{equation}\label{#1}}

\title{On the quantum argument shift method}

\author{Yasushi Ikeda,\quad Alexander Molev \quad and\quad Georgy Sharygin}

\date{} 
\maketitle


\begin{abstract}
In a recent work by two of us
the argument shift method was extended
from the symmetric algebra $\Sr(\g)$ of the general linear Lie algebra $\g$ to the
universal enveloping algebra $\U(\g)$. We show in this paper that
some features of this `quantum argument shift method' can be applied to the remaining
classical matrix Lie algebras $\g$. We prove that a single application of
the quasi-derivation to central elements of $\U(\g)$ yields elements
of the corresponding quantum Mishchenko--Fomenko subalgebra.
We show that generators of this subalgebra can be obtained by iterated application
of the quasi-derivation to generators of the center of $\U(\g)$.

%

\end{abstract}



\section{Introduction}
\label{sec:int}

Let $\g$ be a simple Lie algebra over $\CC$ with basis
elements $Y_1,\dots,Y_l$.
The symmetric algebra $\Sr(\g)$ is equipped with
the {\em Lie--Poisson bracket} which is determined by the condition that its value
$\{Y_i,Y_j\}$ on the basis elements coincides with the commutator $[Y_i,Y_j]$ in $\g$.

Let $P=P(Y_1,\dots,Y_l)$ be an element of $\Sr(\g)$ of a certain degree $d$.
Fix any element $\mu\in\g^*$ and let $t$ be a variable.
Make the substitution $Y_i\mapsto Y_i+t\ts\mu(Y_i)$
and expand as a polynomial in $t$,
\beql{polmf}
P\big(Y_1+t\ts\mu(Y_1),\dots,Y_l+t\ts\mu(Y_l)\big)
=P^{(0)}+P^{(1)} t+\dots+P^{(d)} t^d
\eeq
to define elements $P^{(i)}\in \Sr(\g)$ associated with $P$ and $\mu$.
Denote by
$\overline\Ac_{\mu}$ the subalgebra
of $\Sr(\g)$ generated by all elements $P^{(i)}$
associated with all $\g$-invariants $P\in \Sr(\g)^{\g}$.
The subalgebra $\overline\Ac_{\mu}$ of $\Sr(\g)$ is known
as the {\it Mishchenko--Fomenko subalgebra\/} or
{\it shift of argument subalgebra\/}. Its key property observed in \cite{mf:ee} states
that $\overline\Ac_{\mu}$ is Poisson commutative; that is, $\{R,S\}=0$
for any elements $R,S\in\overline\Ac_{\mu}$.

This property led Vinberg~\cite{v:sc} to raise the {\em quantization problem} asking
whether it is possible to
construct a commutative
subalgebra $\Ac_{\mu}$ of $\U(\g)$ which `quantizes' $\overline\Ac_{\mu}$
in the sense that $\gr\ts\Ac_{\mu}=\overline\Ac_{\mu}$.
In the case where $\mu\in\g^*$ is regular semisimple,
such quantizations were produced by Nazarov and Olshanski~\cite{no:bs}
with the use of the Yangian for $\g=\gl_N$ and the twisted Yangians
associated with the orthogonal and symplectic Lie algebras $\g$.
A different construction of such commutative
subalgebra of $\U(\gl_N)$ was provided by Tarasov~\cite{t:cs}
via a symmetrization map.

A commutative subalgebra $\Ac_{\mu}$ of $\U(\g)$
for arbitrary simple Lie algebra $\g$ and any $\mu\in\g^*$
was constructed with the use of the
associated {\em Feigin--Frenkel center} $\z(\wh\g)$ which
is a certain commutative subalgebra
of $\U\big(t^{-1}\g[t^{-1}]\big)$.
The subalgebra $\Ac_{\mu}$ was shown by Rybnikov~\cite{r:si} and
Feigin, Frenkel and Toledano Laredo~\cite{fft:gm} to solve
Vinberg's quantization problem for all regular $\mu\in\g^*$. Moreover,
it was conjectured in \cite[Conjecture~1]{fft:gm} that the property
$\gr\ts\Ac_{\mu}=\overline\Ac_{\mu}$ holds for all $\mu$. The conjecture
was proved in \cite{fm:qs} and \cite{my:qn} for the Lie algebras $\g$
of types $A$ and $C$.
Explicit generators of the Feigin--Frenkel center $\z(\wh\g)$
and the {\em quantum Mishchenko--Fomenko subalgebra} $\Ac_{\mu}$
have been produced for all simple Lie algebras except for types $E$ and $F$; see
\cite{m:so}, \cite{m:ss}, \cite{my:qn}, \cite{y:sf} and references therein.

Returning to the Mishchenko--Fomenko subalgebra $\overline\Ac_{\mu}$, note that
the coefficients of the polynomial \eqref{polmf} can also be viewed as the images
of the repeated `directional derivatives' $\overline D_{\mu}$ applied to the
polynomial $P\in \Sr(\g)^{\g}$ so that $P^{(k)}=\overline D_{\mu}^{\tss k} P/k!$.
This leads to another version of Vinberg's
quantization problem asking whether the directional derivatives
themselves can be `quantized' to get certain linear operators
\ben
D_{\mu}:\U(\g)\to \U(\g)
\een
with the property that iterative applications of $D_{\mu}$ to elements of the center
$\Zr(\g)$ of $\U(\g)$ yield elements of the quantum Mishchenko--Fomenko subalgebra $\Ac_{\mu}$.
This problem was solved by two of us in \cite{is:as} for the case $\g=\gl_N$.
The solution relies on the {\em quasi-derivations} $\di_{ij}$
acting on $\U(\gl_N)$ as introduced in \cite{gps:bw} and \cite{ms:lr}, while their
further properties were investigated in \cite{i:qo}, \cite{i:sq}, \cite{s:qd} and \cite{s:qdm}.
The operators $\di_{ij}$ on $\U(\gl_N)$ quantize the partial derivations $\di/\di\tss E_{ji}$
on the symmetric algebra $\Sr(\gl_N)$ regarded as the algebra of polynomials
in the $N^2$ variables $E_{ij}$, where $E_{ij}$ denote the standard basis elements of $\gl_N$.
More precisely, the action of $\di_{ij}$
is determined by the properties
\beql{deini}
\di_{ij}(1)=0,\qquad
\di_{ij}\tss E_{kl}=\de_{kj}\de_{il},
\eeq
and the {\em quantum Leibniz rule}
\beql{qlr}
\di_{ij}\tss (fg)=(\di_{ij}\tss f)g+f(\di_{ij}\tss g)-\sum_{k=1}^N (\di_{ik}\tss f)(\di_{kj}\tss g),
\qquad f,g\in \U(\gl_N).
\eeq
The operators $\di_{ij}$ are well-defined; we will include a verification below in
Proposition~\ref{prop:consiA} for completeness; cf. \cite{gps:bw}.

We will think of $\mu\in\gl_N^*$ as the
$N\times N$ matrix $\mu=[\mu_{ij}]$ with $\mu_{ij}=\mu(E_{ij})$.
Combine the quasi-derivations into the matrix $D=[\di_{ij}]$
and set $D_{\mu}=\tr\ts\mu D$. The main result of \cite{is:as} can be stated as
follows.

\bth\label{thm:typeA}
For any element $z\in \Zr(\gl_N)$ and all natural
powers $p$, the elements $D_{\mu}^p\ts z$ belong to
the subalgebra $\Ac_{\mu}$ of $\U(\gl_N)$.
\eth

We will give a proof of Theorem~\ref{thm:typeA} in Sec.~\ref{sec:qdA}
which relies on the tensor matrix techniques and
simplifies
some calculations
in \cite{is:as}. In fact, an alternative version of the quasi-derivations was used therein,
although the two versions are closely related; see
Remark~\ref{rem:difd} below.

It turns out that unlike the action
of the quasi-derivations $\overline D_{\mu}$
on the subalgebra $\overline\Ac_{\mu}$
of the symmetric algebra, the subalgebra $\Ac_{\mu}$ is not preserved by $D_{\mu}$;
see Remark~\ref{rem:notpre}.

Our goal in this paper is to explore the extension of the quantum argument shift method to
the remaining classical Lie
algebras of types $B,C$ and $D$. Define the orthogonal Lie algebras $\oa_N$ with $N=2n+1$ and $N=2n$
and symplectic Lie algebra $\spa_N$ with $N=2n$,
as subalgebras of $\gl_N$ spanned by the elements $F_{i\tss j}$ with $i,j\in\{1,\dots,N\}$,
\beql{fij}
F_{i\tss j}=E_{i\tss j}-E_{j\pr i\pr}\Fand F_{i\tss j}
=E_{i\tss j}-\ve_i\ts\ve_j\ts E_{j\pr i\pr},
\eeq
respectively, for $\oa_N$ and $\spa_N$.
We use the notation $i\pr=N-i+1$, and
in the symplectic case set
$\ve_i=1$ for $i=1,\dots,n$ and
$\ve_i=-1$ for $i=n+1,\dots,2n$.
We will
denote by $\g_N$ any of the Lie algebras $\oa_N$
or $\spa_N$ and define the quasi-derivations by restriction from $\gl_N$.

\bde\label{def:quadef}
The {\em quasi-derivations} $\di_{ij}$ with $i,j\in\{1,\dots,N\}$ act on $\U(\g_N)$ by the rule $\di_{ij}(1)=0$,
\ben
\di_{ij}\tss F_{kl}=\de_{kj}\de_{il}-\de_{ki'}\de_{j'l}
\een
in the orthogonal case, and
\ben
\di_{ij}\tss F_{kl}=\de_{kj}\de_{il}-\ve_k\ts\ve_l\ts \de_{ki'}\de_{j'l}
\een
in the symplectic case, subject to the quantum Leibniz rule \eqref{qlr}
with $f,g\in \U(\g_N)$.
\qed
\ede

Any element $\mu\in\g_N^*$ will be regarded as the
$N\times N$ matrix $\mu=[\mu_{ij}]$ with the entries $\mu_{ij}=\mu(F_{ij})$.
Set $D_{\mu}=\tr\ts\mu D$ with $D=[\di_{ij}]$.

It turns out that a direct analogue of Theorem~\ref{thm:typeA} does not hold
for the Lie algebras of types $B,C$ and $D$; see Remark~\ref{rem:seca} below.
Nonetheless, we show that Proposition~\ref{prop:actionA}
and Theorem~\ref{thm:syman} do extend to the orthogonal and symplectic Lie
algebras; see Proposition~\ref{prop:actionBCD}
and Theorem~\ref{thm:symanBCD}. In particular,
generators of the subalgebra $\Ac_{\mu}$ of $\U(\g_N)$
can be obtained from certain generators of the center $\Zr(\g_N)$ by repeated
applications of the quasi-derivations $D_{\mu}$.

\section{Quasi-derivations in type $A$}
\label{sec:qdA}

We will combine the basis elements $E_{ij}$
of the Lie algebra $\gl_N$ into the matrix $E=[E_{ij}]$.
For each $a\in\{1,\dots,m\}$
introduce the element $E_a$ of the algebra
\beql{tenprka}
\underbrace{\End\CC^{N}\ot\dots\ot\End\CC^{N}}_m{}\ot\U(\gl_N)
\eeq
by
\beql{matnota}
E_a=\sum_{i,j=1}^{N}
1^{\ot(a-1)}\ot e_{ij}\ot 1^{\ot(m-a)}\ot E_{ij},
\eeq
where the $e_{ij}\in \End\CC^{N}$ denote the matrix units.
Similarly, such subscript notation will be used for other matrices
whose entries are elements of $\U(\gl_N)$ or operators on this algebra.
Accordingly, we can regard the matrix $D=[\di_{ij}]$ as the element
\beql{D}
D=\sum_{i,j=1}^N e_{ij}\ot\di_{ij}.
\eeq
The quantum Leibniz rule \eqref{qlr} takes the matrix form
\beql{ql}
D(fg)=(Df)g+f(Dg)-(Df)(Dg),\qquad f,g\in \U(\gl_N).
\eeq
Using notation \eqref{matnota},
in the algebra \eqref{tenprka} with $m=2$ we get
\beql{doeone}
D_1E_2=\sum_{i,j=1}^N e_{ij}\ot1\ot\di_{ij}\ts\sum_{k,l=1}^N 1\ot e_{kl}\ot E_{kl}
=\sum_{i,j=1}^N e_{ij}\ot e_{ji}\ot1=P_{12}\ot 1,
\eeq
where $P_{12}$ is the permutation operator
\beql{p}
P_{12}=\sum_{i,j=1}^N e_{ij}\ot e_{ji}
\eeq
swapping tensor factors in
$\End\CC^{N}\ot\End\CC^{N}$.

\bpr\label{prop:consiA}
The quasi-derivations $\di_{ij}$ are well-defined on $\U(\gl_N)$.
\epr

\bpf
The defining relations of $\U(\gl_N)$ are well-known to admit the matrix form
\beql{deffA}
E_1\ts E_2-E_2\ts E_1=P_{12}\ts E_2-E_2\ts P_{12};
\eeq
see e.g. \cite[Sec.~4.2]{m:so}. By regarding $\U(\gl_N)$ as the quotient of the tensor algebra $\Tr(\gl_N)$
by these relations, note that properties \eqref{deini}
and \eqref{qlr} determine well-defined operators on $\Tr(\gl_N)$.
We then need to verify that the two-sided ideal of
$\Tr(\gl_N)$ generated by the matrix elements of
\ben
E_1\ts E_2-E_2\ts E_1-P_{12}\ts E_2+E_2\ts P_{12}
\een
is preserved by the operators $\di_{ij}$.
Using an additional copy of $\End\CC^{N}$ in \eqref{tenprka}
labelled by $0$ and
applying $D_0$ with the use of \eqref{ql} we get
\begin{multline}
D_0\big(E_1\tss E_2-E_2\tss E_1-P_{12}\tss E_2+E_2\tss P_{12}\big)\\
=
P_{01}E_2+E_1P_{02}-P_{01}P_{02}-P_{02}E_1-E_2P_{01}+P_{02}P_{01}
-P_{12}\tss P_{02}+P_{02}\tss P_{12}.
\non
\end{multline}
This is zero because
\ben
P_{01}P_{02}
=P_{02}P_{12}\Fand
P_{02}P_{01}=P_{12}P_{02},
\een
thus implying the desired property.
\epf

As a first step in proving Theorem~\ref{thm:typeA}, we need the following.

\bpr\label{prop:actionA}
We have $D_{\mu}\ts z\in\Ac_{\mu}$ for any $z\in\Zr(\gl_N)$.
\epr

\bpf
We will use a particular family of generators of $\Zr(\gl_N)$. We start by evaluating the images
of the powers of the matrix $E$ under the action of the quasi-derivation matrix $D$.

\ble\label{lem:doact}
For any $p\geqslant 1$ the expression $D_0 E_1^p$ equals a linear combination
of elements of the form
\ben
E_0^{k}E_1^{l}\Fand E_0^{k}E_1^{l}P_{01}
\een
with $k+l\leqslant p-1$.
\ele

\bpf
We have
\ben
D_0\ts E_1^p=\sum_{s=1}^p (-1)^{s-1} \sum_{p_i}\ts E_1^{p_0}P_{01}E_1^{p_1}P_{01}\dots P_{01}E_1^{p_s}
=\sum_{s=1}^p (-1)^{s-1} \sum_{p_i}\ts E_1^{p_0}E_0^{p_1}E_1^{p_2}E_0^{p_3}\dots P_{01}^{\ve_s},
\een
where $\ve_s=0$ if $s$ is even and $\ve_s=1$ if $s$ is odd. Now use the well-known relation
\beql{yang}
[E_1^r,E_0^s]=\sum_{a=1}^{\min\{r,s\}}(E_0^{a-1} E_1^{r+s-a}P_{01}-E_0^{r+s-a} E_1^{a-1}P_{01}),
\eeq
pointed out in \cite{o:ea} (with the proof credited to N.~Reshetikhin), where it was used for the centralizer
construction of the Yangians. It allows us
to bring the monomials to the required form.
\epf

Returning to the proof of the proposition, note that
any element of $\Zr(\gl_N)$ is a linear combination of terms
\beql{monce}
\tr_{1,\dots,m}^{}\ts E_1^{p_1}\dots E_m^{p_m}
\eeq
for some $m\geqslant 0$ and natural numbers $p_i$, where $\tr_i^{}$ denotes
the partial trace over the $i$-th copy of $\End\CC^{N}$ in \eqref{tenprka}.
By Lemma~\ref{lem:doact},
the expression
$
D^{}_0\ts E_1^{p_1}\dots E_m^{p_m}
$
is a linear combination of monomials
\ben
E_0^{k_1}E_1^{l_1}P_{01}^{\ve_1}E_0^{k_2}E_2^{l_2}P_{02}^{\ve_2}\dots E_0^{k_m}E_m^{l_m}P_{0m}^{\ve_m},
\een
where each $\ve_i$ is $0$ or $1$.
The image of the operator $D_{\mu}=\tr^{}_{0} \ts\mu^{}_0\tss D_0$
applied to the monomial \eqref{monce} is
therefore
a linear combination of the traces
\beql{rtmo}
\tr^{}_{0,\dots,m} \ts\mu^{}_0 E_0^{k_1}E_1^{l_1}P_{01}^{\ve_1}
E_0^{k_2}E_2^{l_2}P_{02}^{\ve_2}\dots E_0^{k_m}E_m^{l_m}P_{0m}^{\ve_m}.
\eeq
Using the relation
\beql{mup}
\tr^{}_i\ts X^{}_i P^{}_{0\tss i}=\tr_i\ts P^{}_{0\tss i}X^{}_{0}=X^{}_{0}\ot 1
\eeq
and calculating the partial traces $\tr^{}_1,\dots, \tr^{}_m$ we represent
the expressions \eqref{rtmo}
as products of elements of the center $\Zr(\gl_N)$ and traces of the form
$\tr\ts\mu E^p$.
Since the center $\Zr(\gl_N)$ is contained in the subalgebra $\Ac_{\mu}$,
it remains to verify that all elements $\tr\ts\mu E^p$ belong to $\Ac_{\mu}$.
By \cite[Theorem~9.3.1]{m:so}, the coefficients
of the polynomials in $z^{-1}$ given by
\ben
\tr\ts(-\di_z+\mu+Ez^{-1})^p\ts 1,\qquad p\geqslant 1,
\een
generate the algebra $\Ac_{\mu}$. The coefficient of $z^{-p+1}$ is a linear combination
of traces $\tr\tss E^r\mu\tss E^s$ with the condition $r+s\leqslant p-1$. However,
they are expressible as $\Zr(\gl_N)$-linear combinations of the traces $\tr\ts\mu\tss E^q$.
More precisely, we have
\ben
\tr\ts E^r\mu\tss E^s=\tr_{1,2}\ts E_1^r\mu_2\ts E_2^sP_{12}=\tr_{1,2}\ts \mu_2\ts E_1^r E_2^sP_{12}.
\een
Now apply \eqref{yang} to write
\ben
E_1^r E_2^s=E_2^s E_1^r +\sum_{a=1}^{\min\{r,s\}}(E_2^{a-1} E_1^{r+s-a}P_{12}-E_2^{r+s-a} E_1^{a-1}P_{12})
\een
so that
$
\tr\ts E^r\mu\tss E^s\equiv\tr\ts \mu\tss E^{r+s}
$
modulo a $\Zr(\gl_N)$-linear combination of the traces $\tr\ts\mu\tss E^q$ with $q<r+s$.
\epf

\bre\label{rem:altarg}
The last step in the proof of Proposition~\ref{prop:actionA} can be carried over
in a different way. To show that all elements
$\tr\ts\mu E^p$ belong to $\Ac_{\mu}$, one can use
the characterization
of the subalgebra $\Ac_{\mu}$ of $\U(\gl_N)$ given in \cite{r:cs}
in the case where $\mu$ is a generic diagonal matrix with distinct entries. Since
the $GL_N$-orbits of such matrices are dense in the space of $N\times N$ matrices,
this will imply the claim for all $\mu$.

Assuming that the
distinct diagonal entries of $\mu$ are $\mu_1,\dots,\mu_N$,
it is enough to show that $\tr\ts\mu E^p$ commutes with all
elements $E_{ii}$ and $T_i\in\U(\gl_N)$ for $i=1,\dots,N$, where
\beql{tidef}
T_i=\sum_{k\ne i}\ts\frac{E_{ik}E_{ki}}{\mu_i-\mu_k}.
\eeq
Using \eqref{deffA}, write
\ben
\bal[]
[E_0,\tr^{}_1\ts \mu_1 E_1^p]&=\sum_{a=1}^p
\tr^{}_1\ts \mu_1 E_1^{a-1}(P_{01}E_1-E_1 P_{01})E_1^{p-a}\\[0.3em]
{}&=\tr^{}_1\ts \mu_1 P_{01}E_1^{p-1}-\tr^{}_1\ts \mu_1 E_1^{p-1} P_{01}=E_0^{p-1}\mu_0-\mu_0 E_0^{p-1}.
\eal
\een
By taking the $(k,i)$ entry we find that
\ben
[E_{ki},\tr\ts \mu E^p]=(\mu_i-\mu_k)\ts (E^{p-1})_{ki}.
\een
This is zero for $i=k$, while
\ben
[T_i,\tr\ts \mu E^p]=\sum_{k\ne i}\ts \big(E_{ik} (E^{p-1})_{ki}-
(E^{p-1})_{ik} E_{ki}\big)=0
\een
which is seen by extending the sum to $k=i$ and noting that
$E_{ii}$ commutes with $E_{ii}^{p-1}$.
\qed
\ere

To complete the proof
of Theorem~\ref{thm:typeA}, we will use the characterization
of the subalgebra $\Ac_{\mu}$ of $\U(\gl_N)$
as pointed out in Remark~\ref{rem:altarg} and will
assume that $\mu$ is a diagonal matrix with
distinct entries $\mu_1,\dots,\mu_N$.
It is enough to show that the image $D_{\mu}^p\ts z$ commutes with all
elements $E_{ii}$ and $T_i$ for $i=1,\dots,N$.
As in \cite{is:as}, we use induction on $p$ with the trivial base $p=0$.
If $E_{ii}x = xE_{ii}$ for some $x\in \U(\gl_N)$, then applying $D_{\mu}$ to both sides we conclude that
\ben
E_{ii}\tss D_{\mu}x+\mu_i\tss x-\mu_i\tss\di_{ii}\tss x
=(D_{\mu}x) E_{ii}+x\mu_i-(\di_{ii}\tss x)\mu_i,
\een
implying $E_{ii}\tss D_{\mu}x=(D_{\mu}x) E_{ii}$. This proves that
$D_{\mu}^p\ts z$ commutes with $E_{ii}$. Furthermore, by applying
$D_{\mu}=\tr^{}_{1}\ts\mu_1 D_1$ to both sides of the relation
$
[T_i,z]=0,
$
we find that
\ben
[D_{\mu}T_i, z]+[T_i, D_{\mu}z]-\tr^{}_{1}\ts\mu_1[D_1T_i, D_1z]=0.
\een
By Proposition~\ref{prop:actionA}, $D_{\mu}z\in\Ac_{\mu}$ so that
this relation implies
\beql{modz}
\tr^{}_{1}\ts\mu_1[D_1T_i, D_1z]=0.
\eeq
As a next step, we take \eqref{modz} as the induction base to derive
a more general relation
\beql{modzk}
\tr^{}_{1}\ts\mu_1[D_1T_i, D_1 D_{\mu}^pz]=0,\qquad p\geqslant 0,
\eeq
by induction on $p$. Set $x=D_{\mu}^pz$ for a fixed $p$ and apply
$D_{\mu}=\tr^{}_{0}\ts\mu_0 D_0$ to both sides of relation \eqref{modzk}
to get
\ben
\tr^{}_{1}\ts\mu_1[D_1D_{\mu}T_i, D_1 x]+\tr^{}_{1}\ts\mu_1[D_1T_i,D_1 D_{\mu}x]
-\tr^{}_{0,1}\ts\mu_0\mu_1[D_0D_1T_i, D_0D_1 x]=0,
\een
where we used the observation that $D_0D_1=D_1D_0$ which is easily
verified on monomials in $\U(\gl_N)$. Note that $D_{\mu}T_i$ is a constant
depending on $\mu$ and so
$D_1D_{\mu}T_i=0$. Hence, it remains to verify that
\beql{modux}
\tr^{}_{0,1}\ts\mu_0\mu_1[D_0D_1T_i, D_0D_1 x]=0.
\eeq
However, the definition \eqref{tidef} implies that
\ben
D_0D_1T_i=\sum_{k\ne i} \frac{1}{\mu_i-\mu_k}(e_{ik}\ot e_{ki}\ot 1+e_{ki}\ot e_{ik}\ot 1).
\een
Therefore, since $\mu=\sum \mu_i e_{ii}$, we have
\ben
\mu_0\mu_1 D_0D_1T_i=D_0D_1T_i\ts\mu_0\mu_1.
\een
Using the cyclic property of trace we then derive
\ben
\tr^{}_{0,1}\ts\mu_0\mu_1\ts D_0D_1T_i\ts D_0D_1 x=
\tr^{}_{0,1}\ts D_0D_1T_i\ts\mu_0\mu_1\ts D_0D_1 x
=\tr^{}_{0,1}\ts\mu_0\mu_1\ts D_0D_1 x\ts D_0D_1T_i,
\een
thus completing the verification of \eqref{modux} and \eqref{modzk}.

As a final step of the proof of Theorem~\ref{thm:typeA}, suppose that the relation
$
[T_i,D_{\mu}^p\ts z]=0
$
holds for a fixed $p\geqslant 0$, and apply $D_{\mu}=\tr^{}_{0}\ts\mu_0 D_0$ to both sides.
Since $D_{\mu}T_i$ is a constant, by taking into account \eqref{modzk}
we conclude that
$
[T_i,D_{\mu}^{p+1}\ts z]=0
$
completing the proof of Theorem~\ref{thm:typeA}.

\medskip

\bre\label{rem:difd}
The arguments of \cite{is:as} rely on different quasi-derivations $\wh\di_{ij}$.
They are determined by the same
properties \eqref{deini}, but the quantum Leibniz rule is changed
to
\beql{qlralt}
\wh\di_{ij}\tss (fg)=(\wh\di_{ij}\tss f)g+f(\wh\di_{ij}\tss g)
+\sum_{k=1}^N (\wh\di_{kj}\tss f)(\wh\di_{ik}\tss g),
\qquad f,g\in \U(\gl_N).
\eeq
In a matrix form, this rule can be written as
\beql{qlalt}
\wh D^{\ts t}(fg)=(  \wh D^{\ts t}f)g+f( \wh D^{\ts t}g)+(\wh D^{\ts t}f)(\wh D^{\ts t}g),
\qquad f,g\in \U(\gl_N),
\eeq
for the matrix $\wh D=[\wh\di_{ij}]$.
The two families of quasi-derivations are related via the automorphism of $\gl_N$
taking $E_{ij}$ to $-E_{ji}$, while $\wh\di_{ij}$ corresponds to $-\di_{ji}$.
Hence, Theorem~\ref{thm:typeA} holds in the same form for the operators $\tr\ts\mu \wh D$
as proved in~\cite{is:as}.
\qed
\ere

We will now show that generators of the subalgebra $\Ac_{\mu}$ can be obtained by repeated
applications of the quasi-derivations $D_{\mu}$ to some generators of the center $\Zr(\gl_N)$.

Consider the natural action of the symmetric group $\Sym_m$ on $(\CC^{N})^{\ot m}$
by permutations of the tensor factors.
We let $H^{(m)}$ and $A^{(m)}$ denote the respective images of the
symmetrizer $h^{(m)}$ and anti-symmetrizer $a^{(m)}$ in the group algebra $\CC[\Sym_m]$.
The elements $h^{(m)}$ and $a^{(m)}$
are the idempotents defined by
\ben
h^{(m)}=\frac{1}{m!}\ts\sum_{s\in\Sym_m} s
\Fand a^{(m)}=\frac{1}{m!}\ts\sum_{s\in\Sym_m} \sgn s\cdot s.
\een
We will identify $H^{(m)}$ and $A^{(m)}$ with the respective elements
$H^{(m)}\ot 1$ and $A^{(m)}\ot 1$ of the algebra \eqref{tenprka}.
Consider the elements
$\vp^{(k)}_m$ and $\psi^{(k)}_m$ of the quantum Mishchenko--Fomenko algebra $\Ac_{\mu}$
constructed in \cite{fm:qs}:
\ben
\vp^{(k)}_{m}=\tr_{1,\dots,m}\ts A^{(m)} \mu^{}_1\dots \mu^{}_k
\ts E^{}_{k+1}\dots E^{}_m
\een
and
\ben
\psi^{(k)}_{m}=\tr_{1,\dots,m}\ts H^{(m)} \mu^{}_1\dots \mu^{}_k
\ts E^{}_{k+1}\dots E^{}_m
\een
Note that each family $\vp^{(0)}_{m}$ and $\psi^{(0)}_{m}$
with $m=1,\dots,N$
is algebraically independent and generates the center $\Z(\gl_N)$; see e.g. \cite[Ch.~5]{m:so}.
The elements $\vp^{(0)}_{m}$ coincide with the images of the basic invariants
of $\Sr(\gl_N)$ under the symmetrization map; see e.g. \cite{my:qn}.
Therefore the following theorem yields another proof of the main result
of \cite{s:qdm} stating that the elements $D^p_{\mu}\tss\vp^{(0)}_{m}$ pairwise commute.

\bth\label{thm:syman}
Each family
\ben
D^p_{\mu}\tss\vp^{(0)}_{m}\Fand D^p_{\mu}\tss\psi^{(0)}_{m}
\een
with $m=1,\dots,N$ and $p=0,\dots,m-1$
generates
the algebra $\Ac_{\mu}$.
Moreover, if $\mu\in\gl_N^*$ is regular, then
each family is algebraically independent.
\eth

\bpf
According to \cite[Corollary~4.6]{fm:qs}, each family of elements $\vp^{(k)}_m$ and $\psi^{(k)}_m$
with $m=1,\dots,N$ and $k=0,\dots,m-1$ possesses the required properties.
Therefore, it will be sufficient to prove that
under the action of the operator $D_{\mu}$, we have
\ben
D_{\mu}:\vp^{(k)}_{m}\mapsto a_1\vp^{(k+1)}_{m}+a_2\vp^{(k+1)}_{m-1}+\dots+a_{m-k}\vp^{(k+1)}_{k+1}
\een
and
\ben
D_{\mu}:\psi^{(k)}_{m}\mapsto b_1\psi^{(k+1)}_{m}+b_2\psi^{(k+1)}_{m-1}+\dots+b_{m-k}\psi^{(k+1)}_{k+1}
\een
for some nonzero constants $a_i$, $b_i$.
By applying the quantum Leibniz rule,
for the image of $\vp^{(k)}_{m}$ we get
\begin{multline}
\tr^{}_0 \ts\mu^{}_0 D^{}_0\ts\tr_{1,\dots,m}\ts A^{(m)} \mu^{}_1\dots \mu^{}_k
\ts E^{}_{k+1}\dots E^{}_m\\[0.4em]
{}=\tr_{0,1,\dots,m}\ts \mu^{}_0\ts A^{(m)} \mu^{}_1\dots \mu^{}_k
\ts\sum_{s=1}^{m-k}\ts\sum_{k+1\leqslant i_1<\dots<i_s\leqslant m}\ts(-1)^{s-1}
\ts E^{}_{k+1}\dots P^{}_{0i_1}\dots P^{}_{0i_s}\dots E^{}_m,
\non
\end{multline}
where $P^{}_{0i_1},\dots, P^{}_{0i_s}$ replace the respective factors
$E^{}_{i_1},\dots, E^{}_{i_s}$. Observe that
\beql{app}
A^{(m)}\ts P^{}_{0i_1}\dots P^{}_{0i_s}=
A^{(m)}\ts P^{}_{i_1i_2}\dots P^{}_{i_1i_s}\ts P^{}_{0i_1}=(-1)^{s-1}\ts A^{(m)}\ts P^{}_{0i_1}
\eeq
so that
the expression simplifies to
\ben
\tr_{0,1,\dots,m}\ts \mu^{}_0\ts A^{(m)} \mu^{}_1\dots \mu^{}_k
\ts\sum_{s=1}^{m-k}\ts\sum_{k+1\leqslant i_1<\dots<i_s\leqslant m}
\ts E^{}_{k+1}\dots P^{}_{0i_1}\dots \wh{\ }\dots E^{}_m,
\een
where the factors $P^{}_{0i_2}\dots P^{}_{0i_s}$ are omitted.
By \eqref{mup} the image of $\vp^{(k)}_{m}$ takes the form
\ben
\tr_{1,\dots,m}\ts A^{(m)} \mu^{}_1\dots \mu^{}_k
\ts\sum_{s=1}^{m-k}
\ts\sum_{k+1\leqslant i_1<\dots<i_s\leqslant m}\ts \mu^{}_{i_1}
\ts E^{}_{k+1}\dots \wh{\ }\dots \wh{\ }\dots E^{}_m.
\een
Finally, we use the cyclic property of trace and conjugate the product by a suitable
permutation $p\in\Sym_m$ by the rule
\begin{multline}
\tr_{1,\dots,m}\ts A^{(m)} X=\sgn P\cdot \tr_{1,\dots,m}\ts A^{(m)} P X
=\sgn P\cdot \tr_{1,\dots,m}\ts A^{(m)} p(X)P\\[0.4em]
=\sgn P\cdot \tr_{1,\dots,m}\ts P  A^{(m)} p(X)=\tr_{1,\dots,m}\ts A^{(m)} p(X),
\label{conj}
\end{multline}
where $P$ is the image of $p$ in the algebra of endomorphisms. Then
the image of $\vp^{(k)}_{m}$ will be written as
\ben
\tr_{1,\dots,m}\ts A^{(m)} \mu^{}_1\dots \mu^{}_k\mu^{}_{k+1}
\ts\sum_{s=1}^{m-k}\binom{m-k}{s}
\ts E^{}_{k+2}\dots E^{}_{m-s+1}.
\een
This is a linear combination of the required form
\ben
a_1\vp^{(k+1)}_{m}+a_2\vp^{(k+1)}_{m-1}+\dots+a_{m-k}\vp^{(k+1)}_{k+1},
\een
because the partial traces of the anti-symmetrizer are found by
\beql{patra}
\tr^{}_{m-s+2,\dots,m}\ts A^{(m)}=\binom{N-m+s-1}{s-1}\binom{m}{s-1}^{-1}\ts A^{(m-s+1)};
\eeq
see e.g. \cite[Ch.~3]{m:so}.
The same argument applies to the elements $\psi^{(k)}_{m}$.
\epf

\bre\label{rem:notpre}
Note that, in general, the algebra $\Ac_{\mu}$ is not preserved
by the operator $D_{\mu}$. To give an example, consider the product
of two elements $\tr\ts \mu E$ and
$\tr\ts E^3$ of the subalgebra $\Ac_{\mu}\subset \U(\gl_N)$
and apply $D_{\mu}$. Working in the algebra \eqref{tenprka}, we get
\ben
D_{\mu}\ts\tr\ts \mu E\ts\tr\ts E^3=\tr^{}_{0,1,2}\ts\mu_0D_0\ts \mu_1E_1E_2^3.
\een
By using the quantum Leibniz rule \eqref{ql}, we find that
modulo elements of $\Ac_{\mu}$, the image equals
\ben
\tr^{}_{0,1,2}\ts\mu_0\mu_1P_{01}\ts(P_{02}E_2^2+E_2P_{02}E_2+E_2^2 P_{02}).
\een
Due to \eqref{deffA}, this coincides with $3\ts\tr\ts\mu^2E^2$ modulo elements of
$\Ac_{\mu}$. However, in general, the element $\tr\ts\mu^2E^2$ does not belong to $\Ac_{\mu}$
which can be checked
by the same calculation as in Remark~\ref{rem:seca} below;
cf. \cite[Example~9.3.2]{m:so}.
\qed
\ere

\section{Quasi-derivations in types $B,C$ and $D$}
\label{sec:qdBCD}

We will combine the elements introduced in \eqref{fij} into the
$N\times N$ matrix
$F=[F_{ij}]$ and use matrix notation as in \eqref{matnota}.
We set $D=[\di_{ij}]$, where $\di_{ij}$ are the quasi-derivations
introduced in Definition~\ref{def:quadef}.
Relation \eqref{doeone} is now replaced by
\beql{dfphi}
D_1F_2=\Phi_{12}\ot 1,\qquad \Phi_{12}=P_{12}-Q_{12},
\eeq
where $P_{12}$ is defined in \eqref{p}, whereas
\ben
Q_{12}=\sum_{i,j=1}^N e_{ij}\ot e_{i'j'}\Fand
Q_{12}=\sum_{i,j=1}^{2n} \ve_i\tss\ve_j\ts e_{ij}\ot e_{i'j'}
\een
in the orthogonal and symplectic case, respectively.
The following is a counterpart of Proposition~\ref{prop:consiA} which we will
verify by a similar argument. We let $\Zr(\g_N)$ denote the center of $\U(\g_N)$.

\bpr\label{prop:actionBCD}
We have $D_{\mu}\ts z\in\Ac_{\mu}$ for any $z\in\Zr(\g_N)$.
\epr

\bpf
The defining relations of $\U(\g_N)$ can be written as
\beql{deff}
F_1\ts F_2-F_2\ts F_1=\Phi_{12}\ts F_2-F_2\ts \Phi_{12};
\eeq
see e.g. \cite[Ch.~5]{m:so}.
We will use the following counterpart of Lemma~\ref{lem:doact}.

\ble\label{lem:doactbcd}
For any $p\geqslant 1$ the expression $D_0 F_1^p$ equals a $\Zr(\g_N)$-linear combination
of elements of the form
\ben
F_0^{k}F_1^{l},\qquad F_0^{k}F_1^{l}P_{01}\Fand F_0^{k}Q_{01}F_0^{l}
\een
with $k+l\leqslant p-1$.
\ele

\bpf
We have
\ben
D_0 F_1^p=\sum_{s=1}^p (-1)^{s-1}
\sum_{p_i}\ts F_1^{p_0}\Phi_{01}F_1^{p_1}\Phi_{01}\dots \Phi_{01}F_1^{p_s}.
\een
Since $P_{12}Q_{12}=Q_{12}P_{12}=\pm Q_{12}$, we can write
\ben
\Phi_{01}F_1^{q}\Phi_{01}=(P_{01}-Q_{01})F_1^{q}(P_{01}-Q_{01})
=F_0^{q}\mp F_0^{q}Q_{01}\mp Q_{01}F_0^{q}+Q_{01} F_1^q\ts Q_{01},
\een
taking the upper sign in the orthogonal case
and the lower sign in the symplectic case.
The operator $Q_{01}$ has the properties
\beql{qxqtr}
Q_{01} F_0^q\ts Q_{01}=Q_{01} F_1^q\ts Q_{01}=Q_{01}\ts\tr\ts F^q
\eeq
and
\beql{qprf}
Q_{01} F_1^q=Q_{01} (F_0^q)',\qquad F_1^q Q_{01} = (F_0^q)'Q_{01},
\eeq
where the prime denotes the matrix transposition associated with the form defining
the orthogonal and symplectic Lie algebras:
\ben
(X\pr)_{ij}=\begin{cases} X_{j'i'}\qquad&\text{in the orthogonal case,}\\
\ve_i\tss\ve_j\ts X_{j'i'}\qquad&\text{in the symplectic case.}
\end{cases}
\een
Now we will use Olshanski's twisted Yangians~\cite{o:ty}
and their relation with the orthogonal and symplectic Lie algebras.
The {\em twisted Yangian} $\Y(\g_N)$ is generated by the entries
of the family of $N\times N$ matrices $S^{(r)}, r=1,2,\dots$ subject to
the defining relations
\ben
\Big(1-\frac{P_{01}}{u-v}\Big)\ts S_0(u)\ts \Big(1+\frac{Q_{01}}{u+v}\Big)\ts S_1(v)=
S_1(v)\ts \Big(1+\frac{Q_{01}}{u+v}\Big)\ts S_0(u)\ts \Big(1-\frac{P_{01}}{u-v}\Big)
\een
and
\beql{symrel}
S(-u)^{\tss\prime}=S(u)\pm\frac{S(u)-S(-u)}{2u},
\eeq
where
\ben
S(u)=1+\sum_{r=1}^{\infty}S^{(r)}u^{-r}.
\een
The {\em evaluation homomorphism} is the epimorphism $\Y(\g_N)\to \U(\g_N)$ defined by
\ben
S(u)\mapsto 1+F\big(u\pm \frac12\big)^{-1}.
\een
Furthermore, the mapping
\ben
S(u)\mapsto c(u)\ts S(-u-N/2)^{-1}
\een
defines an automorphism of the algebra $\Y(\g_N)$, where $c(u)$ is a power series
in $u^{-1}$ with coefficients in the center of $\Y(\g_N)$; see also
\cite[Ch.~2]{m:yc} for a more detailed discussion of these results.
By twisting the evaluation homomorphism with this automorphism,
we get another homomorphism $\Y(\g_N)\to \U(\g_N)$ defined by
\ben
S(u)\mapsto \si(u)\ts \Big(1-F\big(u+\frac{N\pm 1}2\big)^{-1}\Big)^{-1},
\een
where $\si(u)$ is a certain power series
in $u^{-1}$ with coefficients in the center $\Zr(\g_N)$ of $\U(\g_N)$. The images
of the generator matrices then take the form
\ben
S^{(r)}\mapsto F^r+\text{a\ts\  $\Zr(\g_N)$-linear combination of $F^q$ with $q< r$}.
\een
The twisted Yangian defining relations imply that
\beql{prex}
(F^r)'=(-1)^r F^r+\text{a\ts\  $\Zr(\g_N)$-linear combination of $F^q$ with $q<r$},
\eeq
which follows from \eqref{symrel},
whereas
\beql{detw}
F_1^r\ts F_0^s=\text{a\ts\  $\Zr(\g_N)$-linear combination
of $F_0^k F_1^l,\ \  F_0^k F_1^l\ts P_{01}$\ \ and\ \ $F_0^{k}Q_{01}F_0^{l}$},
\eeq
as follows with the use of \eqref{qprf} and \eqref{prex}.
A repeated use of these relations together with \eqref{qxqtr} and
the observation that the traces $\tr\ts F^q$ belong to $\Zr(\g_N)$
complete the proof of the lemma.
\epf

We now continue proving Proposition~\ref{prop:actionBCD} assuming
that $\g_N$ is of type $B$ or $C$. Then
any element of $\Zr(\g_N)$ is a linear combination of terms
\beql{moncebcd}
\tr_{1,\dots,m}^{}\ts F_1^{p_1}\dots F_m^{p_m}
\eeq
for some $m\geqslant 0$ and natural numbers $p_i$. By Lemma~\ref{lem:doactbcd},
the expression
\ben
D^{}_0\ts F_1^{p_1}\dots F_m^{p_m}
\een
is a $\Zr(\g_N)$-linear combination of products
\begin{multline}
\big(a_1F_0^{k_1}F_1^{l_1}+b_1F_0^{k_1}F_1^{l_1}P_{01}+c_1F_0^{k_1}Q_{01}F_0^{l_1}\big)\\[0.4em]
{}\times{}\dots \times\big(a_mF_0^{k_m}F_m^{l_m}+b_mF_0^{k_m}F_m^{l_m}P_{0m}+c_mF_0^{k_m}Q_{0m}F_0^{l_m}\big),
\non
\end{multline}
for some constants $a_i,b_i,c_i$.
The application of the operator $D_{\mu}$ to the monomial \eqref{moncebcd} is
therefore
a linear combination of terms
\begin{multline}
\tr^{}_{0,\dots,m} \ts\mu^{}_0 \big(a_1F_0^{k_1}F_1^{l_1}
+b_1F_0^{k_1}F_1^{l_1}P_{01}+c_1F_0^{k_1}Q_{01}F_0^{l_1}\big)\\[0.4em]
{}\times{}\dots \times\big(a_mF_0^{k_m}F_m^{l_m}+b_mF_0^{k_m}F_m^{l_m}P_{0m}+c_mF_0^{k_m}Q_{0m}F_0^{l_m}\big).
\non
\end{multline}
By calculating the partial traces $\tr^{}_1,\dots, \tr^{}_m$ we represent such terms
as products of elements of the center $\Zr(\g_N)$ and the traces of the form
$\tr\ts\mu F^p$. Since $\Zr(\g_N)$ is contained in $\Ac_{\mu}$,
it remains to verify that $\tr\ts\mu F^p$ belongs to the subalgebra $\Ac_{\mu}$ for all
natural powers $p$. To this end, we will use the argument outlined in
Remark~\ref{rem:altarg} in type $A$.
We will assume that $\mu$
is a diagonal matrix whose entries $\mu_1,\dots,\mu_n$
are nonzero and distinct. By the skew-symmetry of $\mu$ we have
$\mu_{i'}=-\mu_i$
for $i=1,\dots,n$, and also $\mu_{n+1}=0$ for $N=2n+1$.
By the results of \cite{r:cs}, the subalgebra $\Ac_{\mu}$ of $\U(\g_N)$
for such generic $\mu$
consists of the elements which commute with all
$F_{ii}$ and $T_i\in\U(\g_N)$ for $i=1,\dots,n$, where
\beql{tidefbcd}
T_i=\sum_{k\ne i}\ts\frac{F_{ik}F_{ki}}{\mu_i-\mu_k}.
\eeq
Using \eqref{deff}, we derive
\ben
\bal[]
[F_0,\tr^{}_1\ts \mu_1 F_1^p]
=\tr^{}_1\ts \mu_1 \Phi_{01}F_1^{p-1}&-\tr^{}_1\ts \mu_1 F_1^{p-1} \Phi_{01}\\[0.3em]
{}&=F_0^{p-1}\mu_0-\mu_0 F_0^{p-1}+\mu_0 (F_0^{p-1})'-(F_0^{p-1})'\mu_0.
\eal
\een
By taking the $(k,i)$ entry we find that
\ben
[F_{ki},\tr\ts \mu F^p]=(\mu_i-\mu_k)\ts \big((F^{p-1})_{ki}-(F^{p-1})_{i'k'}\big)
\een
in the orthogonal case, and
\ben
[F_{ki},\tr\ts \mu F^p]=(\mu_i-\mu_k)\ts \big((F^{p-1})_{ki}-\ve_i\ve_k (F^{p-1})_{i'k'}\big)
\een
in the symplectic case.
The commutators are zero for $i=k$, and continuing with the orthogonal case we get
\beql{titr}
[T_i,\tr\ts \mu F^p]=\sum_{k\ne i}\ts \big(F_{ik} (F^{p-1})_{ki}-
(F^{p-1})_{ik} F_{ki}-F_{ik} (F^{p-1})_{i'k'}+(F^{p-1})_{k'i'}F_{ki}\big).
\eeq
Now extend the sum to $k=i$ and use
the skew-symmetry property $F_{ik}=-F_{k'i'}$ together with
the well-known commutation
relations between the entries of the matrices $F$ and $F^{p-1}$; see e.g. \cite[Corollary~8.9.4]{m:yc}.
We obtain
\ben
F_{ik} (F^{p-1})_{i'k'}=-(F^{p-1})_{i'k'}F_{k'i'} -
(F^{p-1})_{k'k'}+(F^{p-1})_{i'i'}
\een
and
\ben
(F^{p-1})_{k'i'}F_{ki}=-F_{i'k'}(F^{p-1})_{k'i'} -
(F^{p-1})_{k'k'}+(F^{p-1})_{i'i'},
\een
implying that the sum in \eqref{titr} is zero.
The same calculation in the symplectic case
shows that all elements
$\tr\ts\mu F^p$ belong to $\Ac_{\mu}$, thus
completing the proof
of the proposition for the Lie algebras of types $B$ and $C$.

In the remaining case where $\g_N=\oa_{2n}$ is of type $D$, any element of
$\Zr(\oa_{2n})$ is a linear combination of terms \eqref{moncebcd} and
the traces
\beql{moncepf}
\tr_{1,\dots,m}^{}\ts F_1^{p_1}\dots F_m^{p_m}\ts\Pf\ts F
\eeq
for some $m\geqslant 0$ and natural numbers $p_i$, where $\Pf\ts F\in \Zr(\oa_{2n})$
is the {\em Pfaffian} defined by
\beql{pf}
\Pf\ts F
=\frac{1}{2^nn!}\sum_{\si\in\Sym_{2n}}\sgn\si\cdot
F_{\si(1)\ts\si(2)'}\dots
F_{\si(2n-1)\ts\si(2n)'}.
\eeq

The above argument used in type $B$
applies to the elements \eqref{moncebcd} in type $D$ as well,
so we only need to consider the
images of the elements \eqref{moncepf} under the action
of $D_{\mu}$. It will be more convenient to work with the `canonical' presentation of the Lie algebra
$\oa_{2n}$ so that its elements are skew-symmetric matrices.
We can now think of $F$ as the skew-symmetric matrix
with the entries $F_{ij}=E_{ij}-E_{ji}$. Recall the explicit
isomorphism given in \cite[Sec.~8.1]{m:so}
between the presentation defined by \eqref{fij} and the canonical
presentation, which relates
the associated objects by the formulas
\ben
F\mapsto AFA^{-1},\qquad \mu\mapsto A\mu A^{-1},\qquad D\mapsto ADA^{-1},\qquad
\Pf\ts F\mapsto\det A\cdot\Pf\ts F,
\een
where the matrix $A$ is defined by the relation $AA^t=[\de_{ij'}]$.
Accordingly,
the action of
the quasi-derivations $\di_{ij}$ is changed to
\ben
\di_{ij}\tss F_{kl}=\de_{kj}\de_{il}-\de_{ki}\de_{jl},
\een
while the quantum Leibniz rule \eqref{ql} takes the same form
for the matrix $D=[\di_{ij}]$.
We will now regard $\mu\in \oa_{2n}^*$ as the
matrix $\mu=[\mu_{ij}]$ with $\mu_{ij}=\mu(F_{ij})$.
The Pfaffian in the canonical presentation takes a simpler form
\beql{pfcanon}
\Pf\ts F=\sum_{\si}\sgn\si\cdot
F_{\si(1)\ts\si(2)}\dots F_{\si(2n-1)\ts\si(2n)},
\eeq
summed over the elements $\si$ of the subset of $\Sym_{2n}$
which consists of the permutations with the properties
$\si(2k-1)<\si(2k)$ for all $k=1,\dots,n$ and $\si(1)<\si(3)<\dots<\si(2n-1)$.
All coefficients $\pi^{}_{(k)}$ of the polynomial
\beql{polpf}
\Pf\ts \big(\mu+F\tss z^{-1}\big)=\pi^{}_{(0)}z^{-n}+
\dots+\pi^{}_{(n-1)}z^{-1}
+\pi^{}_{(n)}
\eeq
belong to the algebra $\Ac_{\mu}$.

Introduce the $2n\times 2n$ matrix $\Pi=[\di_{ij} \Pf\ts F]$ with entries in
$\U(\oa_{2n})$. Since the subscripts of the generators $F_{ij}$
occurring in the monomials in \eqref{pfcanon} are all distinct, the
rule \eqref{ql} applied for
the calculation of the entries $\di_{ij} \Pf\ts F$ becomes the usual Leibniz rule
containing only linear terms in the quasi-derivations.

By repeating the above calculation performed for types $B$ and $C$, we find that
modulo elements of the subalgebra $\Ac_{\mu}$,
the image
\ben
\tr^{}_0 \ts\mu^{}_0 D^{}_0\ts\tr_{1,\dots,m}^{}\ts F_1^{p_1}\dots F_m^{p_m}\ts \ts\Pf F
\een
is
a linear combination of terms
\begin{multline}
\tr^{}_{0,\dots,m} \ts\mu^{}_0 \big(a_1F_0^{k_1}F_1^{l_1}
+b_1F_0^{k_1}F_1^{l_1}P_{01}+c_1F_0^{k_1}Q_{01}F_0^{l_1}\big)\\[0.4em]
{}\times{}\dots \times\big(a_mF_0^{k_m}F_m^{l_m}+b_mF_0^{k_m}F_m^{l_m}
P_{0m}+c_mF_0^{k_m}Q_{0m}F_0^{l_m}\big)\ts\Pi_0.
\non
\end{multline}
By calculating the partial traces $\tr^{}_1,\dots, \tr^{}_m$ we represent such terms
as products of elements of the center $\Zr(\oa_{2n})$ and the traces of the form
$\tr\ts\mu F^k\Pi$. It remains to verify that the latter belong
to the subalgebra $\Ac_{\mu}$. Note that $\tr\ts\mu\tss \Pi=-2\pi^{}_{(n-1)}\in\Ac_{\mu}$.
Hence it will be sufficient
to verify the relation
\ben
F\Pi=-\Pf\ts F\cdot 1.
\een
The $(i,j)$ entry of the matrix $F\Pi$ equals
\ben
\sum_{k=1}^{2n}F_{ik}\di_{kj}\ts\Pf\ts F=-\sum_{k=1}^{2n}F_{ik}\di_{jk}\ts\Pf\ts F.
\een
This coincides with $-\Pf\ts F$ for $i=j$, whereas
\ben
\sum_{k=1}^{2n}F_{ik}\di_{jk}\ts\Pf\ts F=0\qquad\text{for $i\ne j$}.
\een
To verify the latter relation, note that for any permutation $\pi\in\Sym_{2n}$ the mapping
$
F_{ij}\mapsto F_{\pi(i)\ts\pi(j)}
$
defines an automorphism of the Lie algebra $\oa_{2n}$. Therefore, it suffices to take
$i=1$ and $j=2$; the relation in this case follows from the Pfaffian formula \eqref{pfcanon}.
\epf

\bre\label{rem:seca}
Theorem~\ref{thm:typeA} does not extend to types $B,C,D$ for the quasi-derivations $D_{\mu}$.
To give an example, consider the element $(\tr\ts F^2)^3\in\Z(\oa_N)$ and show that its image
under $D^2_{\mu}$ does not belong to $\Ac_{\mu}$ for $N\geqslant 5$ and general $\mu$.
For the first application of $D_{\mu}$ we find
\beql{fira}
D_{\mu}\ts (\tr\ts F^2)^3=\tr^{}_{0,1,2,3}\ts \mu_0 D_0 F_1^2F_2^2F_3^2.
\eeq
Note that
\ben
D_0 F_i^2=\Phi_{0i}F_i+F_i\Phi_{0i}-\Phi_{0i}^2
\een
and hence
\ben
\tr^{}_i\ts D_0 F_i^2=4\tss F_0-(2N-2).
\een
Therefore, by replacing the labels of the endomorphism
spaces $\End\CC^{N}$, the image \eqref{fira} can be written as
\begin{multline}
3\tss\tr^{}_{1,2,3}\ts \mu_1\big(4\tss F_1-(2N-2)\big)\tss F_2^2F_3^2\\
{}-3\tss\tr^{}_{1,2}\ts \mu_1\big(4\tss F_1-(2N-2)\big)^2\tss F_2^2
+\tr^{}_{1}\ts \mu_1\big(4\tss F_1-(2N-2)\big)^3.
\non
\end{multline}
Take into account the relation $2\ts\tr\ts \mu F^2 = (N-2)\tss\tr\ts \mu\tss F$ and observe
that the images
\ben
D_{\mu}\ts \tr^{}_{1}\ts \mu_1\tss F_1^3\Fand
D_{\mu}\ts \tr^{}_{1,2}\ts \mu_1\tss F_1F_2^2
\een
belong to $\Ac_{\mu}$. Indeed, checking this for the first image, we get
\ben
\tr^{}_{0,1}\ts \mu_0 D_0 \mu_1F_1^3=\tr^{}_{0,1}\ts \mu_0\mu_1(\Phi_{01}F_1^2+
F_1\Phi_{01}F_1+F_1^2\Phi_{01})+\text{linear terms in $F$}.
\een
The linear terms in $F$ are elements of $\Ac_{\mu}$, while
modulo elements of $\Ac_{\mu}$
the first component simplifies to
\ben
4\ts \tr\ts \mu^2 F^2+ 2\ts \tr\ts (\mu F)^2
\een
which belongs to $\Ac_{\mu}$; see \cite[Example~9.4.2]{m:so}. Hence,
calculating as in the proof of Lemma~\ref{lem:doactbcd}, we find that
modulo elements of $\Ac_{\mu}$,
the image of \eqref{fira} under $D_{\mu}$ divided by $12$ equals
\ben
\tr^{}_{0,1,2,3}\ts \mu_0\mu_1\tss F_1\tss F_2^2F_3^2\equiv 32\ts\tr\ts \mu^2 F^2.
\een
Now we use the characterization
of the subalgebra $\Ac_{\mu}$ of $\U(\oa_N)$ in the case where $\mu$
is a generic regular semisimple element of $\oa_N^*$ \cite{r:cs} and
calculate the commutators of $\tr\ts \mu^2 F^2$ with the elements $T_i$
given by \eqref{tidefbcd}.
We will assume that $\mu$
is a diagonal matrix whose entries $\mu_1,\dots,\mu_n$
are nonzero and distinct. Using the defining relations \eqref{deff},
we get
\ben
[F_0,\tr^{}_1\ts \mu_1^2 F_1^2]=
\tr^{}_1\ts \mu_1^2 \big((\Phi_{01}F_1-F_1\Phi_{01})F_1+F_1(\Phi_{01}F_1-F_1\Phi_{01})\big)
\een
which equals
\ben
2\tss F_0^2\tss\mu_0^2-2\tss\mu_0^2\tss F_0^2+(N-2)\tss(\mu_0^2\tss F_0-F_0\tss \mu_0^2).
\een
Hence, by taking the $(k,i)$ entry we find that
\ben
[F_{ki},\tr^{}_1\ts \mu_1^2 F_1^2]=(\mu_i^2-\mu_k^2)\big(2\tss (F^2)_{ki}-(N-2) F_{ki}\big).
\een
Therefore,
\ben
[T_i,\tr\ts \mu^2 F^2]=\sum_{k\ne i}\ts 2\tss(\mu_i+\mu_k)\tss
\big(F_{ik}(F^2)_{ki}-(F^2)_{ik}F_{ki}\big)
\een
which is nonzero for $N\geqslant 5$. For instance, for $N=5$ and $i=1$ the coefficient of
$2\tss\mu_2$ equals
\ben
F_{12}F_{23}F_{31}-F_{13}F_{32}F_{21}+F_{12'}F_{32}F_{31}-F_{13}F_{23}F_{2'1}\ne 0.
\een
\qed
\ere

We will now show that generators of the subalgebra $\Ac_{\mu}$ can be obtained by repeated
applications of the quasi-derivations $D_{\mu}$ to some generators of the center $\Zr(\g_N)$.
We will regard $\mu\in \g_N^*$ as the
$N\times N$ matrix $\mu=[\mu_{ij}]$ with $\mu_{ij}=\mu(F_{ij})$. Introduce the notation
\ben
\ga_m(\om)=\frac{\om+m-2}{\om+2\tss m-2}.
\een
According to \cite[Ch.~9]{m:so} and \cite{my:qn}, the elements
\beql{vpkp}
\vp^{(k)}_{m}=\ga_{m}(\om)\ts\tr_{1,\dots,m}\ts S^{(m)} \mu^{}_1\dots \mu^{}_k
\ts F^{}_{k+1}\dots F^{}_m
\eeq
belong to the quantum Mishchenko--Fomenko algebra $\Ac_{\mu}$, where $S^{(m)}$
is the symmetrizer arising from the action of the Brauer algebra
on the tensor space $(\CC^{N})^{\ot m}$
and is given by
\ben
S^{(m)}=\frac{1}{m!}
\prod_{1\leqslant a<b\leqslant m}
\Big(1+\frac{P_{ab}}{b-a}-\frac{Q_{ab}}
{N/2+b-a-1}\Big),
\een
in the orthogonal case, and by
\ben
S^{(m)}=\frac{1}{m!}
\prod_{1\leqslant a<b\leqslant m}
\Big(1-\frac{P_{ab}}{b-a}+\frac{Q_{ab}}
{n-b+a+1}\Big),
\een
in the symplectic case,
where the products are taken in the lexicographic order
on the pairs $(a\ts b)$ and
we assume the respective specialisations
$\om=N$ and $\om=-2n$.
The factor $\ga_{m}(\om)$ in \eqref{vpkp} is essential only in the symplectic case
for the values $n+1\leqslant m\leqslant 2n$, where
an additional justification is required to make sense of this formula; see \cite[Ch.~5]{m:so}.

In all cases, the elements $\vp^{(0)}_{m}$ belong to the center $\Z(\g_N)$.
For types $B$ and $C$, the family $\vp^{(0)}_{m}$
with $m=2,4,\dots,2n$
is algebraically independent and generates the center, while in type $D$ the elements
$\vp^{(0)}_{m}$
with $m=2,4,\dots,2n-2$ together with the Pfaffian $\Pf\ts F$ defined in \eqref{pf}
are algebraically independent and generate the center \cite[Ch.~5]{m:so}.

\bth\label{thm:symanBCD}
The family
\beql{dpphi}
D^p_{\mu}\tss\vp^{(0)}_{m},\qquad p=0,1,\dots,m-1,
\eeq
with $m=2,4,\dots,2n$
generates
the algebra $\Ac_{\mu}$ in the cases $B$ and $C$, while the family \eqref{dpphi}
with $m=2,4,\dots,2n-2$ together with the elements $D^p_{\mu}\pi^{(0)}$
for $p=0,1,\dots,n-1$
generate
the algebra $\Ac_{\mu}$ in the case $D$.

Moreover, if $\mu\in\g_N^*$ is regular, then
each family is algebraically independent.
\eth

\bpf
The claims will follow from the respective properties of the elements
$\vp^{(p)}_{m}$ and $\pi^{(p)}$. Namely, the family
\beql{dpphiphi}
\vp^{(p)}_{m},\qquad p=0,1,\dots,m-1,
\eeq
with $m=2,4,\dots,2n$
generates
the algebra $\Ac_{\mu}$ in the cases $B$ and $C$, while the family \eqref{dpphiphi}
with $m=2,4,\dots,2n-2$ together with the elements $\pi^{(p)}$
for $p=0,1,\dots,n-1$
generate
the algebra $\Ac_{\mu}$ in the case $D$. The families are
algebraically independent for regular $\mu$; see
\cite[Ch.~9]{m:so} and \cite{my:qn}.

We start by showing that
under the action of the operator $D_{\mu}=\tr\ts\mu D$, we have
\beql{dmua}
D_{\mu}:\vp^{(k)}_{m}\mapsto c^{}_1\vp^{(k+1)}_{m}+c^{}_3\vp^{(k+1)}_{m-2}+\dots
+c^{}_{2r+1}\vp^{(k+1)}_{m-2r}
\eeq
for some constants $c^{}_i$ with $c_1=2(m-k)$, where $r=\lfloor{\frac{m-k-1}{2}}\rfloor$.
We will be assuming first
that $m\leqslant n$ in the symplectic case.
Then the argument for all types follows the same steps as that of Theorem~\ref{thm:syman}
as outlined below, and it relies on
some properties of the symmetrizer
$S^{(m)}$ \cite[Ch.~1]{m:so}. Instead of \eqref{app}, we get the expression
\ben
S^{(m)}\ts (P^{}_{0i_1}-Q^{}_{0i_1})\dots (P^{}_{0i_s}-Q^{}_{0i_s})
=S^{(m)}\ts\Phi_{0i_1}\dots \Phi_{0i_s}
\een
which we deal with as follows. Write
\ben
\Phi_{0i_1}\dots \Phi_{0i_s}=(P^{}_{0i_1}-Q^{}_{0i_1})\Phi_{0i_2}\dots \Phi_{0i_s}
=P^{}_{0i_1}\Phi_{0i_2}\dots \Phi_{0i_s}-Q^{}_{0i_1}\Phi_{0i_2}\dots \Phi_{0i_s}
\een
and observe that
\beql{pqphph}
\begin{aligned}
P^{}_{0i_1}\Phi_{0i_2}\dots \Phi_{0i_s}&=\Phi_{i_1i_2}\dots \Phi_{i_1i_s}P^{}_{0i_1},\\[0.4em]
Q^{}_{0i_1}\Phi_{0i_2}\dots \Phi_{0i_s}&=(-1)^{s-1}Q^{}_{0i_1}\Phi_{i_1i_s}\dots \Phi_{i_1i_2},
\end{aligned}
\eeq
where we used the relations
\ben
Q_{0i}P_{0j}=Q_{0i}Q_{ij}\Fand Q_{0i}Q_{0j}=Q_{0i}P_{ij},\qquad i\ne j.
\een
Hence,
\ben
S^{(m)}\ts P^{}_{0i_1}\Phi_{0i_2}\dots \Phi_{0i_s}=S^{(m)}\ts P^{}_{0i_1}.
\een
Under the trace over all endomorphism algebras, the product $\Phi_{i_1i_s}\dots \Phi_{i_1i_2}$
can be moved in front of $S^{(m)}$ to give
\ben
\Phi_{i_1i_s}\dots \Phi_{i_1i_2}S^{(m)}=S^{(m)}.
\een
Furthermore, in addition to \eqref{mup}, by the skew-symmetry of $\mu$, we also have
\beql{trqmu}
\tr_0\ts \mu^{}_0 Q^{}_{0i_1}=-\tr_0\ts \mu^{}_{i_1}Q^{}_{0i_1}=-1\ot\mu^{}_{i_1}.
\eeq
This yields the required value of $c_1$ in \eqref{dmua}.
In the final step, for $s\geqslant 2$ calculate the partial traces of the symmetrizer
over the endomorphism spaces $m-s+2,\dots,m$ by
using
\ben
\tr_m\ts\ga_m(\om)\ts S^{(m)}=\pm\tss\frac{\om+m-2}{m}\ts
\ga_{m-1}(\om)\ts S^{(m-1)},
\een
with the plus and minus sign taken in the
orthogonal and symplectic case, respectively; see  \cite[Sec.~1.3]{m:so}. The summand
\ben
\tr^{}_{1,\dots,m-s+1}\ts S^{(m-s+1)} \mu^{}_1\dots \mu^{}_{k+1}
\ts F^{}_{k+2}\dots F^{}_{m-s+1}
\een
will occur in the expansion with the coefficient containing the factor $1-(-1)^{s}$
so that only odd values of $s$ make a contribution in \eqref{dmua}.

Now return to the symplectic case, where we give an alternative proof covering
all values of $m$ with $1\leqslant m\leqslant 2n$. By making use of the formula \cite[(3.16)]{my:qn},
we get an equivalent expression for the elements \eqref{vpkp},
\beql{pkla}
\vp^{(k)}_{m}=\tr_{1,\dots,m}\ts A^{(m)} \mu^{}_1\dots \mu^{}_k
\ts F^{}_{k+1}\dots F^{}_m,
\eeq
which is valid
for all $m=1,\dots,2n$. Using the quantum Leibniz rule
together with \eqref{dfphi}, we obtain
\begin{multline}
\tr^{}_0 \ts\mu^{}_0 D^{}_0\ts\tr_{1,\dots,m}\ts A^{(m)} \mu^{}_1\dots \mu^{}_k
\ts F^{}_{k+1}\dots F^{}_m\\[0.4em]
{}=\tr_{0,1,\dots,m}\ts \mu^{}_0\ts A^{(m)} \mu^{}_1\dots \mu^{}_k
\ts\sum_{s=1}^{m-k}\ts\sum_{k+1\leqslant i_1<\dots<i_s\leqslant m}\ts(-1)^{s-1}
\ts F^{}_{k+1}\dots \Phi^{}_{0i_1}\dots \Phi^{}_{0i_s}\dots F^{}_m,
\non
\end{multline}
where $\Phi^{}_{0i_1},\dots, \Phi^{}_{0i_s}$ replace the respective factors
$F^{}_{i_1},\dots, F^{}_{i_s}$.
Now we use relations \eqref{pqphph} and then calculate the trace over the copy of the endomorphism
algebra $\End\CC^{N}$ labelled by $0$ with the use
of \eqref{mup} and \eqref{trqmu} to bring the expression to the form
\begin{multline}
\tr_{1,\dots,m}\ts A^{(m)} \mu^{}_1\dots \mu^{}_k
\ts\sum_{s=1}^{m-k}\ts\sum_{k+1\leqslant i_1<\dots<i_s\leqslant m}\ts(-1)^{s-1}
\ts F^{}_{k+1}\dots\wh{\ }\dots \wh{\ }\dots F^{}_m\Phi^{}_{i_1i_2}\dots \Phi^{}_{i_1i_s}\mu_{i_1}\\
{}+\tr_{1,\dots,m}\ts A^{(m)} \mu^{}_1\dots \mu^{}_k
\ts\sum_{s=1}^{m-k}\ts\sum_{k+1\leqslant i_1<\dots<i_s\leqslant m}
\tss F^{}_{k+1}\dots\wh{\ }\dots \wh{\ }\dots F^{}_m\ts\mu_{i_1}\tss\Phi^{}_{i_1i_s}\dots \Phi^{}_{i_1i_2}.
\non
\end{multline}
Furthermore, for $s\geqslant 2$
\ben
P^{}_{i_1i_s}\mu_{i_1}=\mu_{i_s}P^{}_{i_1i_s}\Fand P^{}_{i_1i_s}A^{(m)}=-A^{(m)},
\een
whereas
\ben
Q^{}_{i_1i_s}\mu_{i_1}=-Q^{}_{i_1i_s}\mu_{i_s}=Q^{}_{i_1i_s}P^{}_{i_1i_s}\mu_{i_s}
=Q^{}_{i_1i_s}\mu_{i_1}P^{}_{i_1i_s},
\een
so that the product $\Phi^{}_{i_1i_s}\mu_{i_1}$ in the first sum can be
replaced by $-\mu_{i_s}$. Similarly, the product $\mu_{i_1}\tss\Phi^{}_{i_1i_s}$ in the second sum
can also be replaced by $-\mu_{i_s}$. Therefore, by applying
the property \eqref{conj} for
suitable permutations $p$ we can write the image of the element
$\vp^{(k)}_{m}$ under the action of the quasi-derivative
$D_{\mu}$ as
\begin{multline}
\tr_{1,\dots,m}\ts A^{(m)} \mu^{}_1\dots \mu^{}_k\tss\mu^{}_{k+1}
\ts\sum_{s=1}^{m-k}\big(1-(-1)^s\big)\binom{m-k}{s}
\ts F^{}_{k+2}\dots F^{}_{m-s+1}\\
{}\times \Phi^{}_{m-s+2,m-s+3}\dots \Phi^{}_{m-s+2,m},
\non
\end{multline}
where the empty products for $s=1$ and $s=2$ are considered to be equal to $1$.
The following lemma is easily verified by induction.

\ble\label{lem:antphi}
For $r\leqslant m$ we have the relations
\begin{align}
\non
A^{(m)}\Phi^{}_{12}\dots \Phi^{}_{1r}&=A^{(m)}(1+Q_{23})(1+Q_{45})\dots (1+Q_{r-1,r})\quad
&&\text{\rm if $r$ is odd, and}\\[0.4em]
A^{(m)}\Phi^{}_{12}\dots \Phi^{}_{1r}&=-A^{(m)}(1+Q_{23})(1+Q_{45})\dots (1+Q_{r-2,r-1})(1+Q_{1,r})\quad
&&\text{\rm if $r$ is even.}
\non
\end{align}
\ele

Using the lemma and working with odd values of $s$ we come to evaluating
the traces of the form
\beql{treval}
\tr_{1,\dots,m}\ts A^{(m)} \mu^{}_1\dots \mu^{}_k\tss\mu^{}_{k+1}
\ts F^{}_{k+2}\dots F^{}_{m-s+1}
Q^{}_{p,p+1}Q^{}_{p+2,p+3}\dots Q^{}_{p+2q,p+2q+1},
\eeq
where $p=m-s+2$ and $p+2q+1\leqslant m$. We will show that
this trace equals
\ben
\tr_{1,\dots,m}\ts A^{(m-s+1)} \mu^{}_1\dots \mu^{}_k\tss\mu^{}_{k+1}
\ts F^{}_{k+2}\dots F^{}_{m-s+1},
\een
up to a constant factor. If $p+2q+1<m$, then we will calculate the partial traces
of the anti-symmetrizer $A^{(m)}$ over the copies of the endomorphism
algebra $\End\CC^{N}$ labelled by the symbols
$p+2q+2,\dots,m$ with the use of \eqref{patra}.
Therefore, we may assume that $p+2q+1=m$ and argue by induction on $m-p$.
Now use the formula
\ben
A^{(m)}=\frac{1}{m}\big(1-P_{1m}-\dots-P_{m-1,m}\big)\tss A^{(m-1)}
\een
and the cyclic property of trace to write \eqref{treval}
as a linear combination of the terms
\ben
\tr_{1,\dots,m}\ts A^{(m-1)} \mu^{}_1\dots \mu^{}_k\tss\mu^{}_{k+1}
\ts F^{}_{k+2}\dots F^{}_{m-s+1}
Q^{}_{p,p+1}\dots Q^{}_{m-1,m}\tss P^{}_{i\tss m},
\een
with $i=1,\dots,m$ assuming $P_{m\tss m}=1$. Since $Q^{}_{m-1,m}\tss P^{}_{m-1, m}=-Q^{}_{m-1,m}$
and $\tr_m \ts Q^{}_{m-1,m}=1$, the cases $i=m-1$ and $i=m$ are taken care of by the induction hypothesis.
For $i\leqslant m-2$, we have
\ben
\tr_m \ts Q^{}_{m-1,m}\tss P^{}_{i\tss m}=\tr_m \ts P^{}_{i\tss m}\tss Q^{}_{i,m-1}=Q^{}_{i,m-1}.
\een
If $i\leqslant k+1$, then
\ben
A^{(m-1)} \mu_i\tss Q^{}_{i,m-1}=-A^{(m-1)} \mu_{m-1}\tss Q^{}_{i,m-1}=
A^{(m-1)} \mu_{m-1}\tss P^{}_{i,m-1}\tss Q^{}_{i,m-1}=-A^{(m-1)} \mu_i\tss Q^{}_{i,m-1}
\een
showing that the corresponding term is zero. The same property takes place for
the values of $i$ with $k+1<i\leqslant m-s+1$.
Finally, if $m-s+1<i\leqslant m-2$, then use one of the relations
\ben
Q^{}_{i,i+1}\tss Q^{}_{i,m-1}=Q^{}_{i,i+1}\tss P^{}_{i+1,m-1}\qquad\text{or}\qquad
Q^{}_{i-1,i}\tss Q^{}_{i,m-1}=Q^{}_{i,i+1}\tss P^{}_{i-1,m-1}.
\een
It remains to use the cyclic property of trace to eliminate the factor $P^{}_{i-1,m-1}$
and apply the induction hypothesis.
\epf

This completes
the proof of the theorem in the symplectic case. It remains to consider
the even orthogonal case and the elements arising from the Pfaffian.
As we pointed out above,
the rule \eqref{ql} applied for
the calculation of the entries $\di_{ij} \Pf\ts F$ in the canonical presentation \eqref{pfcanon},
contains only linear terms in the quasi-derivations. Therefore,
the elements $D^p_{\mu}\pi^{(0)}$ coincide with the respective coefficients
$\pi^{(p)}$ in the expansion \eqref{polpf}, up to nonzero multiplicative factors.
This completes the proof of the theorem for type $D$ as well.

\bigskip\bigskip

\small

\noindent
Y.I. \& G.S.:\newline
Faculty of Mechanics and Mathematics\\
Lomonosov Moscow State University\\
GSP-1, Leninskie Gory, 119991, Moscow, Russia

\vspace{3 mm}

\noindent
G.S.:\newline
Kurchatov Institute \\
1 Kurchatov Square, 123182, Moscow, Russia

\vspace{3 mm}

\noindent
A.M.:\newline
School of Mathematics and Statistics\newline
University of Sydney,
NSW 2006, Australia\newline
alexander.molev@sydney.edu.au

\end{document}